\documentclass[12pt]{amsart}
\usepackage{amscd,amssymb}
\usepackage{graphicx}
\usepackage[arrow,matrix]{xy}
\usepackage{pst-all}

\theoremstyle{plain}
\newtheorem{thm}[subsection]{Theorem}

\newtheorem{prop}[subsection]{Proposition}
\newtheorem{cor}[subsection]{Corollary}
\newtheorem{con}[subsection]{Conjecture}

\theoremstyle{definition}
\newtheorem{rk}[subsection]{Remark}

\newtheorem{ex}[subsection]{Example}

\numberwithin{equation}{section} 

\newcommand{\Z}{\mathbb{Z}}
\newcommand{\Q}{\mathbb{Q}}

\newcommand{\C}{\mathbb{C}}

\newcommand{\PP}{\mathbb{P}}

 \DeclareMathOperator{\rank}{rank}
\DeclareMathOperator{\im}{Im}

 \DeclareMathOperator{\mult}{mult}

\DeclareMathOperator{\codim}{codim}

\DeclareMathOperator{\Rad}{Rad}

\begin{document}

\title[Polar Cremona Transformations and Monodromy of Polynomials]{Polar Cremona Transformations and Monodromy of Polynomials}
\author[Imran Ahmed]{Imran Ahmed}

\address{School of Mathematical Sciences,
         Government College University,
         68-B New Muslim Town Lahore,
         PAKISTAN.}
\email {iahmedthegr8@gmail.com}



\begin{abstract}
Consider the gradient map associated to any non-constant homogeneous
polynomial $f\in \C[x_0,\ldots,x_n]$ of degree $d$, defined by
\[\phi_f=grad(f)\colon D(f)\to \PP^n, (x_0\colon \ldots \colon
x_n) \to (f_0(x)\colon \ldots \colon f_n(x))\]where $D(f)=\{x\in
\PP^n; f(x)\neq 0\}$ is the principal open set associated to $f$ and
$f_i=\frac{\partial f}{\partial x_i}$. This map corresponds to polar
Cremona transformations. In Proposition \ref{p1}
 we give a new lower
bound for the degree $d(f)$ of $\phi_f$ under the assumption that the
projective hypersurface $V:f=0 $ has only isolated
singularities. When $d(f)=1$, Theorem \ref{t4} yields very strong conditions on the singularities of $V$.
\end{abstract}

\maketitle

%\tableofcontents

\section{Introduction}

Consider the gradient map associated to any non-constant homogeneous
polynomial $f\in\C[x_0,\ldots,x_n]$ of degree $d$, defined by
\[\phi_f=grad(f)\colon D(f)\to \PP^n, (x_0\colon \ldots \colon
x_n) \to (f_0(x)\colon \ldots \colon f_n(x))\]where $D(f)=\{x\in
\PP^n; f(x)\neq 0\}$ is the principal open set associated to $f$ and
$f_i=\frac{\partial f}{\partial x_i}$. This map corresponds to the polar
Cremona transformations considered by Dolgachev in \cite{Do}, see also
\cite{DP}, \cite{D5}, \cite{CRS}, \cite{FP}.

In section 2 we recall basic facts on  the degree $d(f)$
of the gradient map $grad(f)$, emphasizing in Theorem \ref{t2}
the relation to the Bouquet Theorem of L\^e in \cite{Le}.

In section 3 we consider the challenging Conjecture \ref{con} due to A. Dimca and S. Papadima in \cite{DP},
saying roughly that $d(f)>1$ for most projective hypersurface
$V: f=0$ having only isolated singularities.
In Proposition \ref{p1} we find a new lower bound for the degree of
the gradient map $\phi_f$ in this case.

In the last section we show how the monodromy of a polynomial function $h:\C^n \to \C$,
naturally associated to $f$, gives extremely strong conditions on the monodromy operators associated to the singularities of $V$
if $d(f)=1$.

\section{The Degree of The Gradient}

Let $d(f)=deg(\phi_f)$ denote the degree of the gradient map, which is defined as follows. For a
dominant map $\phi_f$ we have the following equivalent
definitions:

\medskip

\noindent (i) There is a Zariski open and dense subset $U$ in
$\PP^n$ such that for all $u\in U$ the fiber $\phi_f^{-1}(u)$ has
exactly $d(f)$ points;

\medskip

\noindent (ii) the rational fraction field extension
$\phi_{f}^{*}\colon K(\PP^n)\to K(D(f))$ has degree $d(f)$, see
Mumford \cite{M}, Proposition (3.17).

In particular, this latter formulation implies that $d(f)=1$ if and only if the
gradient map $\phi_f$ induces a birational isomorphism of the
projective space $\PP^n$.\\
The degree of the gradient map $\phi_f$ is defined to be zero if the
gradient map is not dominant.\\
Note that in all the above we may replace the open set $D(f)$ by the
larger open set
\[E(f)=\PP^n\backslash \{x\in\PP^n;\,\,f_0(x)=f_1(x)=\ldots=f_n(x)=0\}\]
without changing the degree of the gradient map (to see this just
see description (i) given above for the degree).\\
One has the following topological description of the degree
$d(f)$ of the gradient map $grad(f)$, see \cite{DP}.

\begin{thm}\label{t1}
For any non-constant homogeneous polynomial $f\in
\C[x_0,\ldots,x_n]$, the complement $D(f)$ is homotopy equivalent to
a CW complex obtained from $D(f)\cap H$ by attaching $d(f)$ cells of
dimension $n$, where $H$ is a generic hyperplane in $\PP^n$. In
particular, we have
\[d(f)=(-1)^n\chi (D(f)\backslash H)\]
\end{thm}

Note that the meaning of 'generic' here is quite explicit: the
hyperplane $H$ has to be transversal to a stratification of the
projective hypersurface $V$.\\
This yields in particular the following corollary, see \cite{DP}  and also \cite{FP} for a recent, completely different approach.

\begin{cor}\label{c1}
The degree of the gradient map $grad(f)$ depends only on the reduced
polynomial $f_{red}$ associated to $f$.
\end{cor}

Moreover, Theorem \ref{t1} can be restated in the following way,
which shows that for any projective hypersurface $V$, if we choose
the hyperplane at infinity $H$ in a generic way, then the topology
of the affine part $X=V\backslash H$ is very simple. For details, see \cite{DP}.

\begin{thm}\label{t2}
For any non-constant homogeneous polynomial $f\in
\C[x_0,\ldots,x_n]$, the affine part $X(f)=V(f)\backslash H$ of the
corresponding projective hypersurface $V(f)$ with respect to a
generic choice of the hyperplane at infinity $H$ is homotopy
equivalent to a bouquet of $(n-1)$-spheres. The number of spheres in
this bouquet is the degree $d(f)$.
\end{thm}

\begin{rk}\label{r1}
 (i) Using Thom's second Isotopy Lemma, see for instance \cite{D2}, it
follows that, for any projective variety $V$, the topology of affine
part $X=V\backslash H$ is independent of $H$, for a generic
hyperplane $H$. For this reason, we will use the alternative simpler
notation $V_a$ for the generic affine piece $X$ of the projective
variety $V$. Exactly the same argument as in the proof of Theorem
\ref{t2} shows that $V_a$ is homotopy equivalent to a bouquet of
$k$-spheres
when $V$ is a complete intersection of dimension $k$.\\
(ii) It is not difficult to construct projective hypersurfaces $V$ with
isolated singularities such that for a given hyperplane $H_0$,
$V\backslash H_0$ is smooth and contractible, see \cite{CD2}. However,
since $H_0$ is not a generic hyperplane in this case, this does not
imply $d(f)=0$.
\end{rk}

\section{Hypersurfaces with Isolated Singularities}

In this section we consider the following conjecture, see \cite{DP}, end of section 3, and \cite{D5}.
\begin{con}\label{con}
Let $f\in \C[x_0,\ldots,x_n]$ be a reduced homogeneous polynomial
such that\\
(a) $d=deg(f)>2$ and $n>2$;\\
(b) the associated projective hypersurface $V(f)$ has only isolated
singularities.\\
Then $d(f)\neq 1$.
\end{con}

The following result was obtained by A. Dimca in \cite{D5}.

\begin{thm}\label{t3}
This conjecture is true if either

\noindent (i) all the singularities of the hypersurface
$V(f)$ are  weighted homogeneous, or

\noindent (ii) the hypersurface
$V(f)$ is a $\Q$-manifold.
\end{thm}

Using Theorem \ref{t2} and known facts on the topology of special
fibers in a deformation of an isolated hypersurface singularity, we
have
\begin{equation}\label{eq1}
d(f)=(d-1)^n-\mu(V(f)),
\end{equation}
where $\mu(V(f))$ is the sum of the Milnor numbers of all the
singularities of $V(f)$, see \cite{D2}, p. 161 for details. When all
these singularities are weighted homogeneous, then
$\mu(V(f))=\tau(V(f))$, where $\tau(V(f))$ is the sum of the Tjurina
numbers of all the singularities of $V(f)$. The claim (i) above follows from deep results by
du Plessis and Wall \cite{dPW} giving upper bounds for $\tau(V(f))$.

The proof of the second claim in Theorem \ref{t3} is much easier, and is generalized
in the proof of Proposition \ref{p1} below.

The following example shows that the conjecture is optimal.
\begin{ex}\label{e1}
An example with $d(f)=2$ can be obtained as follows.
Let $n=3$, $d=3$ and let $f$ be the equation of a cubic surface with
singularities $A_1A_5$ or $E_6$ (for the existence of cubic surfaces having these configurations of singularities, see
\cite{BW}).
Now, using equation  \eqref{eq1}, we get $d(f)=(3-1)^3-6=2$ in either
case.
\end{ex}

\begin{prop}\label{p1}
Let $f\in \C[x_0,\ldots,x_n]$ be a homogeneous polynomial of degree
$d>2$ and such that $n\geq 3$. If the associated projective
hypersurface $V=V(f)\subset \PP^n$ has only isolated singularities,
say at the points $a_1,\ldots,a_p$, then
$$d(f)\geq b_{n-2}^0(W_{n-2}^d)-\mu^0(V).$$
Here $b_{n-2}^0(W_{n-2}^d)$ is the primitive middle Betti number of
a smooth projective hypersurface $W_{n-2}^d$ of degree $d$ and
dimension $n-2$ and\[\mu^0(V)=\sum_{j=1,p}\mu^0(V,a_j)\]is the sum
of the ranks of the radicals of the Milnor lattices
$L_j$ corresponding to the singularities $(V,a_j)$.
\end{prop}

Note that the hypersurface
$V=V(f)$ is a $\Q$-manifold if and only if $\mu^0(V)=0$ and that $b_{n-2}^0(W_{n-2}^d)>1$.
In this way Proposition \ref{p1} implies the second claim in Theorem \ref{t3}.

\begin{proof}
Let $L$ be the Milnor lattice of the isolated singularity obtained as the vertex of the affine cone over a generic hyperplane section $W=V \cap H_0$ of the hypersurface $V$.
 Consider the
lattice morphism
$$\phi_V:\,L_1\oplus L_2\oplus\ldots\oplus L_p\stackrel{\psi_V}{\hookrightarrow}L\stackrel{p}{\to}\overline{L}=L/\Rad L$$
as in \cite{D2} p.161, where $\psi_V$ is the obvious (primitive) embedding of lattices,
$$\Rad L=\{x\in L;\,xy=0,\,\text{for all}\, y\in L\}$$
and $p$ is the canonical projection. It follows that
\begin{equation}\label{eq2}
\ker \phi_V=\Rad L\cap (L_1\oplus L_2\oplus\ldots\oplus L_p).
\end{equation}
Next we show that
\begin{equation}\label{eq3}
\ker\phi_V\subset \Rad\,L_1\oplus \Rad\,L_2\oplus\ldots\oplus \Rad\,L_p.
\end{equation}

Let $v=(v_1,\ldots,v_p)\in \ker\phi_V$. Then, by equation \eqref{eq2},
$v \in \Rad\,L$. By definition, this implies that
$v.w=0$ for all $w\in L$. Taking
$w=(0,\ldots,0,b_j,0,\ldots,0)$ for any $b_j\in L_j$, we get
$v_jb_j=0$ for all $b_j\in L_j$. It implies that $a_j\in \Rad L_j$
for each $j=1,\ldots,p$. It follows that \eqref{eq3} holds.

\medskip

Consider the Milnor number of the singularity $(V,a_k)$ given by
$\mu(V,a_k)=\rank L_k$ and note that  $\mu^0(V)$ denotes the sum
$\sum_{k=1,p}\rank(\Rad\,L_k)$.  Now, by the basic properties of the
rank of a $\Z$-linear mapping between free $\Z$-modules of finite type, we have
$$\rank \im \phi_V=\rank (\oplus_{k=1}^pL_k)-\rank
(\ker\phi_V)=\sum_{k=1,p}\rank\,(L_k)-\rank\,(\ker\phi_V)=$$
$$=\mu(V)-\rank\,(\ker\phi_V).$$
By equation \eqref{eq3}, we get $\rank\,(\ker \phi_V)\leq \mu^0(V)$. Therefore,
we get
\begin{equation}\label{eq4}
\rank \im \phi_V \geq \mu(V)-\mu^0(V)\geq 0.
\end{equation}
Note that $\im \phi_V\subset \overline{L}$, which yields
\begin{equation}\label{eq5}
\rank\overline{L}\geq \rank\im\phi_V\geq\mu(V)-\mu^0(V)\geq 0.
\end{equation}
On the other hand, we have
$$(d-1)^n=b_{n-1}(W_{n-1}^d)+b_{n-2}(W_{n-2}^d)-1,$$
see formula (5.3.27) on p.159 in \cite{D2}. By the definition of the primitive Betti
numbers, this is equivalent to
$$(d-1)^n=b_{n-1}^0(W_{n-1}^d)+b_{n-2}^0(W_{n-2}^d).$$
Using the equations \ref{eq1} and \ref{eq5}, it follows that
$$d(f)=(d-1)^n-\mu(V)\geq b_{n-1}^0(W_{n-1}^d)+b_{n-2}^0(W_{n-2}^d)-(rk\,\overline{L}+\mu^0(V)).$$
Moreover, one has $\rank\overline{L}= b_{n-1}^0(W_{n-1}^d)$, as follows from
Prop (5.3.24), p.157 in \cite{D2}, which implies the required
result.
\end{proof}

For $n=3$, $W_{n-2}^d$ is a smooth projective plane curve of degree
$d$. By the degree-genus formula, we have
$$g=\frac{(d-1)(d-2)}{2}.$$

Thus, we get
$$b_1^0(W_1^d)=(d-1)(d-2).$$

This implies the following.

\begin{cor}\label{cor3.7}
The Conjecture \ref{con} is true for any surface $V$ of degree
$d$ such that $\mu^0(V)<(d-1)(d-2)-1.$
\end{cor}

\section{The use of monodromy of polynomial functions}

In this section we prove the main result, which is a generalization of Proposition \ref{p1}
above, in the key case $d(f)=1$. In order to state it we need some preliminaries.
Assume as above that the hypersurface $V=V(f)\subset \PP^n$ has only isolated singularities,
say at the points $a_1,\ldots,a_p$. At each singular point $a_j$, there is a local Milnor fiber $F_j$
and a monodromy operator $T_j: H_{n-1}(F_j,\C) \to H_{n-1}(F_j,\C)$. Let
\begin{equation}\label{eq6}
\Delta_j(t)= \det (t \cdot I -T_j)
\end{equation}
be the corresponding characteristic polynomial and set
\begin{equation}\label{eq7}
\Delta_V(t)= \prod_{j=1,p}\Delta_j(t).
\end{equation}
For any $\lambda \in \C$, we let $\mult_V (\lambda)$ denote the multiplicity of $\lambda$ as a root of the equation
$\Delta_V(t)=0$. For instance, one clearly has (e.g. by using the proof of Prop. 3.4.7
in \cite{D2}, p.93)
\begin{equation}\label{eq6.5}
\mult_V (1)=\mu^0(V).
\end{equation}

Consider the following basic example.
\begin{ex}\label{e6}
The homogeneous polynomial $g=x_1^d+ \ldots +x_n^d$ has an isolated singularity at the origin
$0 \in \C^n$. Then it is well known that the corresponding characteristic polynomial of the monodromy $T_0$ is given by
$$\Delta_0(t)= \left( \frac{(t^d-1)^{\chi(G)/d}}{t-1}\right)^{(-1)^{n-1}}$$
where $G: g(x)=1$ is the affine Milnor fiber of $g$, and $\chi(G)=1+(-1)^{n-1}\mu(g)$, with
$\mu(g)=(d-1)^n$. Note that $\chi(F)/d=\chi(U)$, where $U=\PP^{n-1} \setminus V(g)$.
Since $V(g)$ is smooth, it follows that $\chi(F)/d=1+(-1)^{n-1}b_{n-2}^0(W_{n-2}^d)$.
\end{ex}
For any $\lambda \in \C$, we let $\mult_0 (\lambda)$ denote the multiplicity of $\lambda$ as a root of the equation
$\Delta_0(t)=0$. It follows that
\begin{equation}\label{eq8}
\mult_0 (\lambda)=b_{n-2}^0(W_{n-2}^d)+(-1)^{n-1}
\end{equation}
when $\lambda$ is a $d$-root of unity, $\lambda \ne 1$, and
\begin{equation}\label{eq9}
\mult_0 (\lambda)=b_{n-2}^0(W_{n-2}^d)
\end{equation}
when $\lambda = 1$.
Now we can state our main result.
\begin{thm}\label{t4}
If $d(f)=1$, then for any $d$-root of unity  $\lambda$ one has
$$\mult_V( \lambda) \geq \mult_0( \lambda)-1.$$
\end{thm}
Using the equality \eqref{eq6.5}, it follows that, by taking $\lambda = 1$,  Theorem \ref{t4} implies
Proposition \ref{p1} when $d(f)=1$.

\proof

Assume that the hyperplane at infinity $H_0:x=0$ in $\PP^n$ is transversal to $V$. Then, in the affine space
$\C^n=\PP^n \setminus H_0$, the corresponding affine part $V_a$ is defined by an equation
$$h(x_1, \ldots,x_n)=f(1, x_1, \ldots,x_n)=0.$$
Since $W=V \cap H_0$ is smooth, it follows that the polynomial $h$ is tame, see  \cite{D2}, p.22. In particular, $h$ has only isolated singularities on $\C^n$ and
$$\sum_x\mu(h,x)=(d-1)^n.$$
If $d(f)=1$, it follows that the polynomial $h$ has precisely $p+1$ isolated singularities:
$p$ of them on the zero fiber $H_1=V_a=h^{-1}(0)$ and the last one, an $A_1$ singularity, on a different fiber, say
$H_2=h^{-1}(b)$, for $b \ne 0$.

The generic fiber $F_h$ of $h$ is diffeomorphic to the fiber $F$ in Example \ref{e6}
(where $g=0$ is an equation for $W$!), and the monodromy at infinity
$T_{\infty}$ of the polynomial $h$ corresponds, under the identification $E=H_{n-1}(F_h,\C)=H_{n-1}(F,\C)$,
to the monodromy operator $T_0$ in Example \ref{e6}.

On the other hand, we have a relation of the type $T_{\infty}=T_1 \circ T_2$, where $T_j$ denotes the monodromy of $h$
about the singular fiber $H_j$, for $j=1,2.$
Let $H=\ker(T_2-I) \subset E$. It is known that $\codim H=1$, see \cite{D3}, p.202, Prop. 6.3.19 (iii).
The monodromy operator $T_{\infty}$ is semisimple, so we get a direct sum decomposition
$$E =\oplus_{\lambda} E_{\lambda}$$
where $E_{\lambda}$ is the $\lambda$-eigenspace corresponding to $T_{\infty}$.

Consider first the case ${\lambda} \ne 1$. Then one clearly has
$$\dim (E_{\lambda} \cap H) \geq \dim (E_{\lambda} )-1=\mult_0({\lambda})-1.$$
On the other hand, for $v \in  E_{\lambda} \cap H$ one has $\lambda v=T_1(v)$, i.e.
$E_{\lambda} \cap H$ is contained in the  eigenspace of  $T_1$ corresponding to the eigenvalue $\lambda$.
It is known that the characteristic polynomial of the monodromy operator $T_1$ is given by
\begin{equation}\label{eq10}
\det(t\cdot I-T_1)=\Delta_V(t) \cdot (t-1)
\end{equation}
see \cite{D3}, p.188, Cor. 6.2.17.
This yields the claim in this case, since obviously $\mult_V({\lambda})$ is greater than the dimension of
eigenspace of  $T_1$ corresponding to the eigenvalue $\lambda$.

In the case ${\lambda} = 1$, it follows from \cite{DN}, Theorem 3.1, (ii), that both restrictions $T_1|E_1$ and  $T_2|E_1$
are the identity of $E_1$. In particular, in this case we have
$$\dim (E_{1} \cap H)=\dim (E_{1} )=\mult_0(1).$$
Due to the factor $(t-1)$ in the formula \eqref{eq10}, this is exactly what is needed to get the claim in this case.

\endproof

\end{document}